\documentclass[12pt]{amsart}
\usepackage{latexsym,amsfonts,amsmath,amssymb,amsthm,url}
\usepackage[english]{babel}
\usepackage[latin1]{inputenc}
\usepackage{amsbsy}
\usepackage{amscd}
\usepackage{graphicx,psfrag,epsfig}
\usepackage{enumerate}
\usepackage{dsfont}
\newtheorem{theorem}{Theorem}
\newtheorem{lemma}[theorem]{Lemma}
\newtheorem{proposition}[theorem]{Proposition}

\newtheorem{remark}[theorem]{Remark}
\usepackage{color}
\newcommand{\RR}{\mathbb R}

\newcommand{\e}{\varepsilon}
\newcommand{\dd}{\;{\rm d}}
\begin{document}

\title[Finite mass self-similar blowing-up solutions]{Finite mass self-similar blowing-up solutions of a chemotaxis system with non-linear diffusion}

\author[A. Blanchet]{Adrien Blanchet$^1$}
\address{$^1$ GREMAQ, CNRS UMR~5604, INRA UMR~1291, Universit\'e de Toulouse, 21 All\'ee de Brienne, F--31000 Toulouse, France}
\email{Adrien.Blanchet@univ-tlse1.fr}

\author[Ph. Lauren\c cot]{Philippe Lauren\c cot$^2$}
\address{$^2$ Institut de Math\'ematiques de Toulouse, CNRS UMR~5219, Universit\'e de Toulouse, F--31062 Toulouse cedex 9, France}
\email{laurenco@math.univ-tlse.fr}

\date{\today}

\thanks{Partially supported by the project EVaMEF ANR-09-JCJC-0096-01.}

\begin{abstract} For a specific choice of the diffusion, the parabolic-elliptic Patlak-Keller-Segel system with non-linear diffusion (also referred to as the quasi-linear Smoluchowski-Poisson equation) exhibits an interesting threshold phenomenon: there is a critical mass $M_c>0$ such that all the solutions with initial data of mass smaller or equal to $M_c$ exist globally while the solution blows up in finite time for a large class of initial data with mass greater than $M_c$. Unlike in space dimension $2$, finite mass self-similar blowing-up solutions are shown to exist in space dimension $d\ge 3$.
\end{abstract}


\maketitle

\section{Introduction}

In space dimension $d=2$, the parabolic-elliptic Patlak-Keller-Segel (PKS) system is a simplified model which describes the collective motion of cells in the following situation: cells diffuse in space and emit a chemical signal, the chemo-attractant, which results in the cells attracting each other. If $\rho$ denotes the density of cells and $c$ the concentration of the chemo-attractant, the PKS system reads \cite{KS71,P53}
\begin{equation}
  \label{PKS}
\left\{
  \begin{array}{l}
\partial_t \rho(t,x)= {\rm div} \left[ \nabla \rho(t,x) - \rho(t,x)\nabla c(t,x)\right]\,, \vspace{.3cm}\\
c(t,x) = (E_2\star\rho)(t,x)\,, \quad E_2(x) = -\displaystyle{\frac{1}{2\pi}\ \ln{|x|}}\,,
  \end{array}
\quad (t,x) \in [0,\infty) \times \RR^2\,.
\right.
\end{equation}
This model may be seen as an elementary brick to understand the aggregation of cells in mathematical biology as it exhibits the following interesting and biologically relevant feature: there is a critical mass above which the density of cells is expected to concentrate near isolated points after a finite time, a property which is related to the formation of fruiting bodies in the slime mold \textit{Dictyostelium discoideum}. Such a phenomenon does not take place if the density of cells is too low. More precisely, given a non-negative integrable initial condition $\rho_0$ with finite second moment, the system \eqref{PKS} has a unique maximal classical solution $(\rho,c)$ defined on some maximal time interval $[0,T)$, $T\in (0,\infty]$. Its first component $\rho$ is non-negative and the mass of $\rho$ (that is, its $L^1$-norm) remains constant through time evolution 
$$
\|\rho(t)\|_1= M :=\|\rho_0\|_1\,, \quad t\in [0,T)\,.
$$
It is well-known that, if $M< 8\,\pi$, the solution to \eqref{PKS} exists globally in time while it blows up in finite time if $M>8\,\pi$, see \cite{BDP06,DP04,Ho03,JL92} and the references therein. More recently, it was shown that there is global existence as well for the critical mass $M=8\,\pi$, the blowup occurring in infinite time with a profile being a Dirac mass of mass $8\pi$ \cite{BCM08}. When the mass $M$ is above $8\pi$, the shape of the finite time blowup  is not self-similar according to asymptotic expansions computed in \cite{CS06,Luxx} (see also \cite{HV96} for a related problem in a bounded domain). In addition, there is no integrable and radially symmetric blowing-up self-similar solution to \eqref{PKS} \cite[Theorem~8]{NS08}.

\smallskip

In space dimension $d \ge 3$, the system~\eqref{PKS} seems to be less relevant from the biological point of view as blowup may occur whatever the value of $M$ \cite{HMV97,Na00}. This means that the diffusion is too weak to balance the aggregation resulting from the chemotactic term. It is however well-known that one can enhance the effect of diffusion to prevent crowding by considering a diffusion of porous medium type which increases the diffusion of the cells when their density $\rho$ is large. This is the generalised version of the Patlak-Keller-Segel model considered in, e.g., \cite{BCL09,CS04,Su06,Su07,ST09}:
\begin{equation}
  \label{nPKS}
\left\{
  \begin{array}{l}
\partial_t \rho(t,x)={\rm div} \left( \nabla \left[\rho^m(t,x)\right] -  \rho(t,x)\nabla c(t,x)\right)\,,\vspace{.3cm}\\
c(t,x) = (E_d\star\rho)(t,x)\,, \quad E_d(x) = c_d\ |x|^{2-d}\,,
  \end{array}
\quad (t,x) \in [0,\infty) \times \RR^d\,,
\right.
\end{equation}
where $m>1$, $c_d:= 1/((d-2)\ \sigma_d)$, and $\sigma_d:=2\ \pi^{d/2}/\Gamma(d/2)$ denotes the surface area of the sphere $\mathbb{S}^{d-1}$ of $\RR^d$. The system \eqref{nPKS} also arises in astrophysics \cite{CS04}  (being then referred to as the generalised Smoluchowski-Poisson equation), and $\rho$ and $c$ denote the density of particles and the gravitational potential, respectively.

For~\eqref{nPKS}, it turns out that there is only one critical exponent of the non-linear diffusion, namely $m_d:=2(d-1)/d$, such that the mass plays a similar role to that in~\eqref{PKS}. Indeed, if $m>m_d$ the diffusion enhancement is too strong and the solutions always exist globally in time whereas if $m<m_d$ the diffusion is not strong enough to compensate the aggregation term and there are solutions blowing up in finite time whatever the value of the mass~\cite{Su06,Su07}. The relevant diffusion is thus achieved in the case when $m=m_d$. In this case, it was proved in~\cite{BCL09} that there is a unique threshold mass $M_c>0$ with the following properties: if the mass $M=\|\rho_0\|_1$ of the initial condition $\rho_0$ is less or equal to $M_c$,  then the corresponding solution to \eqref{nPKS} exists globally in time, whereas given any $M>M_c$ there are initial data $\rho_0$ with mass $M$ such that the corresponding solution blows up in finite time. Thus, for the peculiar choice $m=m_d$ and $d\ge 3$, the system \eqref{nPKS} exhibits the same qualitative behaviour as the PKS system \eqref{PKS} in space dimension $2$. Still, there is a fundamental difference as the latter has no fast-decaying stationary solution with mass $8\pi$ while the former has a two-parameter family of non-negative, integrable, and compactly supported stationary solutions with mass $M_c$ for each $d\ge 3$~\cite[Section~3]{BCL09}.

\smallskip 

It is then tempting to figure out whether this striking difference extends above the critical mass $M_c$ and this leads us to investigate the existence of blowing-up (or backward) self-similar solutions with finite mass. More precisely, since mass remains unchanged throughout time evolution, we look for solutions $(\rho,c)$ to \eqref{nPKS} with $m=m_d$ and $d\ge 3$ of the form
\begin{equation}
\label{spirou}
\rho(t,x)=\frac{1}{s(t)^{d}}\ \Phi\left(\frac{x}{s(t)}\right) \;\;\mbox{ and }\;\; c(t,x)=\frac{1}{s(t)^{d-2}}\ \Psi\left(\frac{x}{s(t)}\right) 
\end{equation}
for $(t,x)\in [0,T)\times\RR^d$ with $s(t):=\left[ d(T-t)\right]^{1/d}$, the time $T$ being an arbitrary positive real number. Note that $s(t)$ converges to zero as $t$ increases to the blowup time $T$.

\smallskip

Our main result is then the following:
\begin{theorem}[Existence of finite mass self-similar blowing-up solutions]\label{th:main}
There exists $M_2 \in (M_c,\infty)$ such that, for any $M$ in $(M_c,M_2]$, there exists at least a non-negative self-similar blowing-up solution $(\rho_M,c_M)$ to~\eqref{nPKS} of the form \eqref{spirou} with a radially symmetric, compactly supported, and non-increasing profile $\Phi_M$ satisfying $\|\rho_M(t)\|_1=\|\Phi_M\|_1=M$ for $t\in [0,T)$ and $\|\rho_M(t)\|_\infty\to\infty$ as $t\to T$.
\end{theorem}

As a consequence of Theorem~\ref{th:main}, we realize that non-negative, integrable, and radially symmetric self-similar blowing-up solutions to \eqref{nPKS} with a non-increasing profile only exist below a threshold mass.  Another by-product of our analysis is the existence of non-negative and non-integrable self-similar blowing-up solutions to~\eqref{nPKS}, see Proposition~\ref{pr:as} below. 

\section{Blowing-up self-similar profiles}\label{bussp}
From now on,
$$
d\ge 3 \;\; \mbox{ and }\;\; m=m_d=\frac{2 (d-1)}{d}\,,
$$
and we look for a solution $(\rho,c)$ to \eqref{nPKS} of the form
\begin{equation}
\label{bu1}
\rho(t,x)=\frac{1}{s(t)^{d}}\ \Phi\left(\frac{x}{s(t)}\right) \;\;\mbox{ and }\;\; c(t,x)=\frac{1}{s(t)^{d-2}}\ \Psi\left(\frac{x}{s(t)}\right) 
\end{equation}
with $s(t)=\left[ d(T-t)\right]^{1/d}$ and $(t,x)\in [0,T)\times\RR^d$ for some given $T>0$. We further assume that $\Phi$ enjoys the following properties:
\begin{equation*}
\left\{
\begin{array}{l}
\Phi\in\mathcal{C}(\RR^d)\cap L^1(\RR^d) \;\mbox{ is radially symmetric and non-negative, }\;\; \\
\\
\Phi^{m-1}\in W^{1,\infty}(\RR^d)\,.
\end{array}
\right.
\end{equation*}
Inserting the ansatz \eqref{bu1} in \eqref{nPKS}  gives that $(\Phi,\Psi)$ solves
\begin{equation*}
\left\{
  \begin{array}{l}
{\rm div} \left( \nabla \left[\Phi^m(y)\right] -  \Phi(y) \nabla \Psi(y) - \Phi(y)\ y\right) = 0\,,\vspace{.3cm}\\
\Psi(y) = (E_d\star\Phi)(y)\,,
  \end{array}
\right.
\end{equation*}
for $y\in\RR^d$. Since $\Psi=E_d \star \Phi$, the radial symmetry of $\Phi$ ensures that of $\Psi$ and, introducing the profiles $(\varphi,\psi)$ of $(\Phi,\Psi)$
\begin{equation}
\label{bu21}
\Phi(y) = \varphi(|y|)\,, \qquad \Psi(y)=\psi(|y|)\,, \qquad y\in\RR^d\,.
\end{equation}
By \cite[Theorem~9.7, Formula~(5)]{LL01}, we have
\begin{equation}
\label{bu3}
\psi(r) = \frac{1}{(d-2) r^{d-2}}\ \int_0^r \varphi(s)\ s^{d-1}\ \dd s + \frac{1}{d-2}\ \int_r^\infty \varphi(s)\ s\ \dd s
\end{equation}
for $r\ge 0$. We can also write the equation for $\varphi$ as
\begin{equation}
\label{bu31}
\partial_r\left( r^{d-1}\ \varphi(r)\ \partial_r J(r) \right) = 0 \;\;\mbox{ with }\;\; J(r):= \frac{2(d-1)}{d-2}\ \varphi^{(d-2)/d}(r) -  \psi(r) - \frac{r^2}{2} \end{equation}
for $r\in (0,\infty)$. Since we are looking for an integrable profile, we formally conclude that 
\begin{equation}
\label{bu4}
\partial_r J(r) = 0 \;\;\mbox{ for }\;\; r\in \mathcal{P}_\varphi:=\{ s\in (0,\infty) \ : \ \varphi(s)>0 \}\,.
\end{equation}
In particular, $J$ is constant on any connected component of $\mathcal{P}_\varphi$. But, if $\mathcal{C}$ is a connected component of $\mathcal{P}_\varphi$, we have either
\begin{equation}
\label{cc1}
\mathcal{C}=(0,R_s) \;\;\mbox{ for some }\;\; R_s\in (0,\infty]\,,
\end{equation}
or
\begin{equation*}
\mathcal{C}= (R_i,R_s) \;\;\mbox{ for some }\;\; R_i\in (0,\infty) \;\;\mbox{ and }\;\; R_s\in (0,\infty]\,. 
\end{equation*}

\begin{remark}
If we additionally assume that the profile $\varphi$ is non-increasing then $\mathcal{P}_\varphi$ has only one connected component which is necessarily of the form \eqref{cc1}.
\end{remark}

Now, take a connected component $\mathcal{C}$ of $\mathcal{P}_\varphi$. It follows from \eqref{bu4} that there is $\mu\in\RR$ such that  
\begin{equation}
\label{bu5}
J(r) = \frac{2(d-1)}{d-2}\ \varphi^{(d-2)/d}(r) -  \psi(r) - \frac{r^2}{2} = - \mu \;\;\mbox{ for }\;\; r\in \mathcal{C}\,.
\end{equation}
Owing to the assumed integrability of $\Phi$, the function $r\mapsto r^{d-1} \varphi(r)$ belongs to $L^1(0,\infty)$ and it follows from \eqref{bu3} that the function $r\mapsto r^{d-2} \psi(r)$ is bounded in $\mathcal{C}$. Therefore \eqref{bu5} only complies with the integrability of $\Phi$ if $R_s<\infty$ which implies the boundedness of $\mathcal{C}$. Introducing
\begin{equation*}
\Xi:=\varphi^{(d-2)/d}
\end{equation*}
 and taking the Laplacian of both sides of \eqref{bu5} yield that $\Xi$ is a positive solution to
\begin{equation}
\label{bu61}
- \frac{d^2 \Xi}{dr^2}(r) - \frac{d-1}{r}\ \frac{d\Xi}{dr}(r) = \frac{d-2}{2(d-1)}\ \left( \Xi(r)^{d/(d-2)} - d \right) \;\;\mbox{ in }\;\; \mathcal{C}\,, 
\end{equation}
with either 
\begin{equation}
\label{bu62}
\partial_r \Xi(0)=\Xi(R_s)=0 \;\;\;\mbox{ if }\;\;\; \mathcal{C}=(0,R_s)
\end{equation} 
or  
\begin{equation}
\label{bu63}
\Xi(R_i)=\Xi(R_s)=0 \;\;\;\mbox{ if }\;\;\; \mathcal{C}=(R_i,R_s)\,.
\end{equation}
 A final change of scale, namely
\begin{equation*}
\eta(r) := \frac{1}{\lambda_d}\ \Xi\left( \frac{r}{\mu_d} \right)\,, \quad \lambda_d:=d^{(d-2)/d}\,, \quad \mu_d:= d^{1/d}\ \left( \frac{d-2}{2(d-1)} \right)^{1/2}\,,
\end{equation*}
leads us to the following boundary-value problem for $\eta$: either
\begin{equation}
\label{bu8}
\left\{
  \begin{array}{l}
\displaystyle{\frac{d^2 \eta}{dr^2}(r) + \frac{d-1}{r}\ \frac{d\eta}{dr}(r) + \eta(r)^{d/(d-2)} - 1 = 0}\,, \quad r\in (0,\mu_d R_s)\,,\vspace{.3cm}\\
\displaystyle{\frac{d\eta}{dr}(0)=0}\,, \quad \eta(\mu_d R_s)=0\,,
  \end{array}
\right.
\end{equation}
or
\begin{equation}
\label{bu9}
\left\{
  \begin{array}{l}
\displaystyle{\frac{d^2 \eta}{dr^2}(r) + \frac{d-1}{r}\ \frac{d\eta}{dr}(r) + \eta(r)^{d/(d-2)} - 1 = 0}\,, \quad r\in (\mu_d R_i,\mu_d R_s)\,,\vspace{.3cm}\\
\displaystyle{\eta(\mu_d R_i)=0}\,, \quad \eta(\mu_d R_s)=0\,.
  \end{array}
\right.
\end{equation}
We have thus reduced our study to one or several boundary-value problems (depending on the number of connected components of $\mathcal{P}_\varphi$) for a nonlinear second order differential equation. The purpose of the next section is then a precise study of this ordinary differential equation. 

However, before going on, let us point out that \eqref{bu61} is not equivalent to \eqref{bu5}. Indeed, since
$$
\partial_r J(r)=\frac{2(d-1)}{d-2}\ \partial_r\Xi(r) + \frac{1}{r^{d-1}}\ \int_0^r \Xi(s)^{d/(d-2)}\ s^{d-1}\ \dd s - r\,, \quad r\in \mathcal{C}\,,
$$
by \eqref{bu31}, the fact that $\Xi$ is a solution to \eqref{bu61} only guarantees that $\partial_r(r^{d-1} \partial_r J(r))=0$ for $r\in\mathcal{C}$. Consequently, there are constants $C_1$ and $C_2$ such that 
$$
\partial_r J(r)=-\frac{(d-2)\ C_1}{r^{d-1}}\,, \quad J(r) = \frac{C_1}{r^{d-2}} + C_2\,, \quad r\in\mathcal{C}\,,
$$
from which \eqref{bu5} follows only if $C_1=0$. On the one hand, if $\mathcal{C}=(R_i,R_s)$ with $0<R_i<R_s$, it is yet unclear whether the boundary conditions \eqref{bu63} might imply this property. On the other hand, if $\mathcal{C}=(0,R_s)$, the boundary conditions \eqref{bu62} ensure that $\partial_r J(0)=0$ and thus $C_1=0$. We shall only deal with this case in the remaining of this paper and thus focus on the non-increasing profiles $\varphi$.

\section{An auxiliary ordinary differential equation}\label{aaode}

For $a\in\RR$, let $u(.,a)\in\mathcal{C}^1([0,r_{\max}(a)))$ denote the maximal solution to the Cauchy problem
\begin{equation}
\label{c1}
\left\{
\begin{array}{l}
\displaystyle{u''(r,a) + \frac{d-1}{r}\ u'(r,a) + |u(r,a)|^{p-1}\ u(r,a) - 1 = 0}\,, \quad r\in [0,r_{\max}(a))\,,\\
\\
u(0,a)= a\,, \quad u'(0,a)=0\,,
\end{array}
\right.
\end{equation}
with $r_{\max}(a)\in (0,\infty]$ and $p=d/(d-2)$.

Clearly, if $a=1$ then $u(.,1)\equiv 1$ is a stationary solution and $r_{\max}(1)=\infty$. We first show that $u(.,a)$ is global for all $a\in\RR$ and oscillates around the value $1$ if $a\ne 1$.

\begin{lemma}\label{lec1}
For each $a\in\RR\setminus\{1\}$, $r_{\max}(a)=\infty$, and the solution $u(.,a)$ to~\eqref{c1} is an oscillatory function in $(0,\infty)$. More precisely,
\begin{itemize}
\item if $a>1$, there is an increasing sequence $(r_i(a))_{i\ge 0}$ of real numbers such that $r_0(a)=0$,
\begin{equation*}
\left\{
\begin{array}{l}
u'(r_i(a),a)=0\,, \quad (-1)^i u'(r,a)<0 \;\;\mbox{ for }\;\; r\in (r_i(a),r_{i+1}(a))\,,\\
\\
u(r_{2i}(a),a)>u(r_{2i+2}(a),a)>1>u(r_{2i+3}(a),a)>u(r_{2i+1}(a),a)
\end{array}
\right.
  \end{equation*}
for $i\ge 0$,
\item if $a<1$, there is an increasing sequence $(r_i(a))_{i\ge 1}$ of real numbers such that $r_1(a)=0$ 
\begin{equation*}
\left\{
\begin{array}{l}
u'(r_i(a),a)=0\,, \quad (-1)^i u'(r,a)<0 \;\;\mbox{ for }\;\; r\in (r_i(a),r_{i+1}(a))\,,\\
\\
u(r_{2i}(a),a)>u(r_{2i+2}(a),a)>1>u(r_{2i+1}(a),a)>u(r_{2i-1}(a),a)
\end{array}   
\right.   
  \end{equation*}
for $i\ge 1$.
\end{itemize}
\end{lemma}

\begin{figure}[h!]
 \begin{minipage}[t]{1\linewidth}
\centering\epsfig{figure=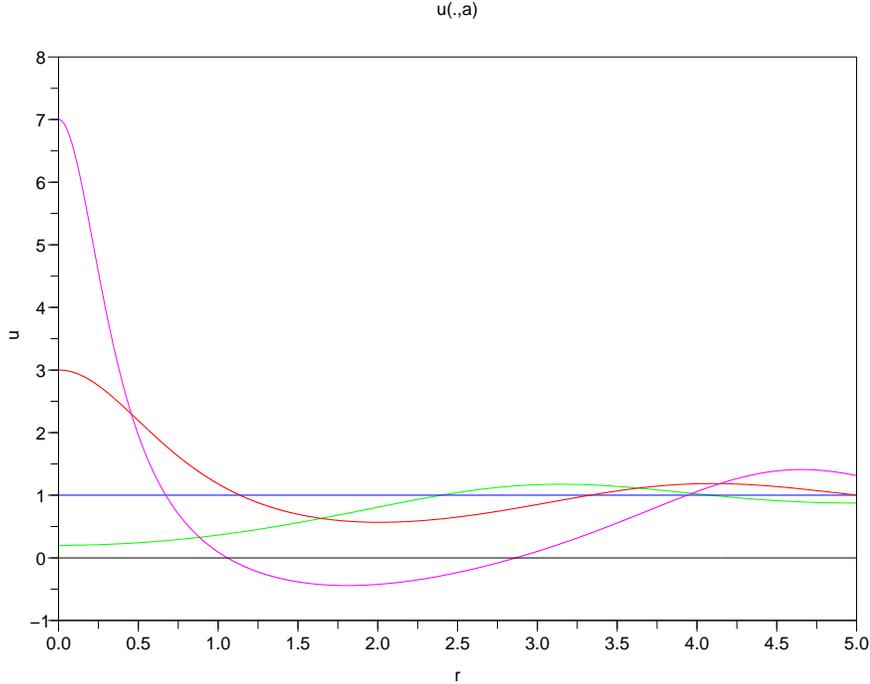,width=10cm,angle=270}
\caption{\small Various oscillating behaviours of $u(.,a)$ for $a\in\{ 0.2 , 1 , 3, 7\}$.}\label{fig:lemma4}
 \end{minipage} \hfill
\end{figure}

These properties are illustrated in Figure~\ref{fig:lemma4}. Notice that, for $a=7$, $u(.,7)$ vanishes at a finite $r$ and thus provides a solution to \eqref{bu8}.  

\begin{proof}[Proof of Lemma~\ref{lec1}] For any $r\in [0,r_{\max}(a))$ consider the functional
\begin{equation}
\label{c2}
E(r,a) := \frac{|u'(r,a)|^2}{2} + \frac{|u(r,a)|^{p+1}}{p+1} - u(r,a) \,.
\end{equation}
By~\eqref{c1}, for all $r\in [0,r_{\max}(a))$
\begin{equation}
\label{c3}
\frac{dE}{dr}(r,a) = - \frac{d-1}{r}\ |u'(r,a)|^2 \le 0\,,
\end{equation}
Obviously $E(r,a) \ge - p/(p+1)$. Owing to \eqref{c3}, $E(r,a)\in [-p/(p+1),E(0,a)]$ for $r\in [0,r_{\max}(a))$ which prevents $u(.,a)$ of becoming unbounded at a finite value of $r$, thereby implying that $r_{\max}(a)=\infty$. We next argue using Sturm's oscillations theorem as in \cite[Lemma~9]{Kw89}, to establish the oscillatory behaviour of $u(.,a)$ for $a\ne 1$.
\end{proof}

\medskip

According to \eqref{bu8}, we are interested in finding solutions to the initial value problem~\eqref{c1} which are positive and vanish at a finite value of $r$. We thus focus on the case $a>0$ and investigate the positivity properties of $u(.,a)$.

\begin{lemma}\label{lec2}
There is a constant $a_c>1$ such that 
\begin{itemize}
\item if $a\in (0,a_c)$, then $u(r,a)>0$ for all $r\ge 0$,
\item if $a=a_c$, then there is $R(a_c)>0$ such that
\begin{equation*}
\left\{
\begin{array}{ll}
u(R(a_c),a_c)=0\vspace{.2cm}\\
u'(R(a_c),a_c)=0\vspace{.2cm}\\
u(r,a_c)>0 \,\quad\mbox{for $r\in [0,R(a_c))$,} 
\end{array}   
\right.    
 \end{equation*}
\item if $a\in (a_c,\infty)$, then there is $R(a)>0$ such that
\begin{equation*}
\left\{
\begin{array}{ll}
u(R(a),a)=0\vspace{.2cm}\\
u'(R(a),a)<0\vspace{.2cm}\\
u(r,a)>0 \,\quad\mbox{for $r\in [0,R(a))$.} 
\end{array}   
\right.   
  \end{equation*}
\end{itemize} 
\end{lemma}

\begin{figure}[h!]
 \begin{minipage}[t]{1\linewidth}
\centering\epsfig{figure=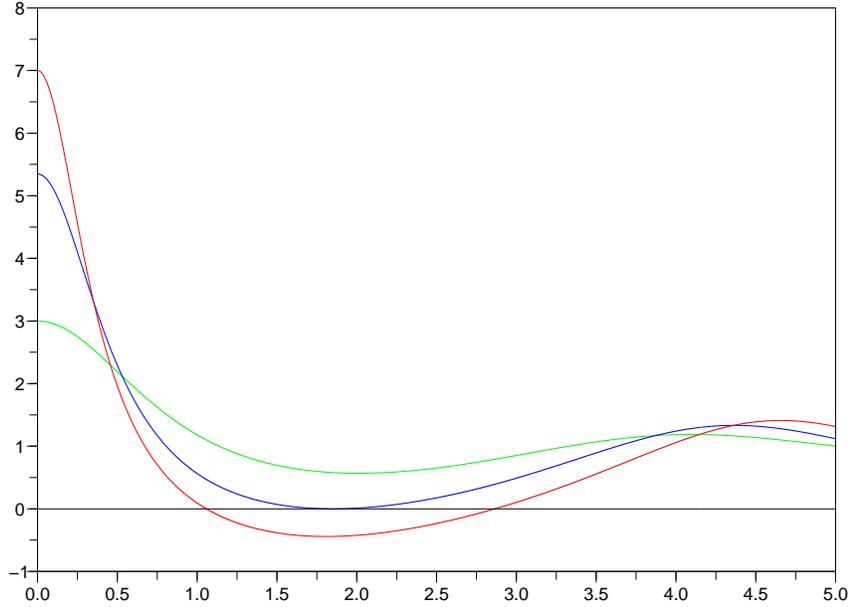,width=10cm,angle=270}
  \caption{\small Behaviour of $u(.,a)$ for $a>a_c$, $a=a_c$ and $a<a_c$.}\label{fig:lemma5}
 \end{minipage} \hfill
\end{figure}

These three possibilities are drawn in Figure~\ref{fig:lemma5}.

\begin{proof}[Proof of Lemma~\ref{lec2}]
For $a>0$, we define 
\begin{equation*}
R(a):= \inf\{ R>0 \ : \ u(r,a)>0 \;\;\mbox{ for }\;\; r\in [0,R) \}\,.
\end{equation*}
Notice that the positivity of $a$ and the continuity of $u(.,a)$ guarantee that $R(a)>0$. We consider the sets
$$
\begin{array}{lcl}
\mathcal{P} & := & \left\{ a>0\ : \ R(a)=\infty \right\}\,, \\
\mathcal{N} & := & \left\{ a>0\ : \ R(a)<\infty \;\mbox{ and }\; u'(R(a),a)<0 \right\}\,, \\
\mathcal{N}_0 & := & \left\{ a>0\ : \ R(a)<\infty \;\mbox{ and }\; u'(R(a),a)=0 \right\}\,.
\end{array}
$$
Clearly, $\mathcal{P} \cup \mathcal{N} \cup \mathcal{N}_0=(0,\infty)$ and $1\in\mathcal{P}$. Actually, if $a\in (0,(p+1)^{1/p})$, then $E(0,a)<0$ and the monotonicity \eqref{c3} of $E$ entails that $E(r,a)<0$ for all $r\ge 0$. But, if $R(a)<\infty$, it readily follows from the definition~\eqref{c2} of the functional $E$ that $E(R(a),a)\ge 0$ whence a contradiction. Therefore, $R(a)=\infty$ for any $a\in (0,(p+1)^{1/p})$ so that 
\begin{equation}
\label{c5}
\left( 0,(p+1)^{1/p} \right) \subset \mathcal{P}\,.
\end{equation}

Consider now $a\in\mathcal{N}_0$. Then $U(x):=u(|x|,a)$ is a radial positive solution to the homogeneous Dirichlet-Neumann free boundary problem $\Delta U + U^p-1=0$ in $B(0,R(a))$ with $U=\partial_\nu U = 0$ on $\partial B(0,R(a))$. According to \cite[Theorem~3~(iii)]{ST00}, there is only one value of $a$ for which this solution has a positive radial solution and it is unique. Consequently, there is a unique $a_c>0$ such that $\mathcal{N}_0=\{a_c\}$.  

Consider next $a\in \mathcal{N} \cup \mathcal{N}_0$ and recall that $a>1$ by \eqref{c5}. Following \cite[Lemma~11]{Kw89}, let us assume for contradiction that there is $\varrho\in (0,R(a))$ such that $u'(\varrho,a)=0$. Either $\varrho\le 1$ and we infer from the definition, the monotonicity of $E$, see~\eqref{c2}-\eqref{c3}, and the definition of $R(a)$ that 
$0>E(\varrho,a)\ge E(R(a),a)\ge 0$ which is a contradiction. Or $\rho>1$ and the oscillating behaviour of the solutions implies, using the notation of~Lemma~\ref{lec1}, that $\varrho\ge r_2(a)$. This implies that $r_1(a)<R(a)$. Then $u(r_1(a),a)\in (0,1)$ and using again \eqref{c2}, \eqref{c3}, and the definition of $R(a)$, we conclude that 
$0>E(r_1(a),a)\ge E(R(a),a)\ge 0$, hence a contradiction. Therefore,
\begin{equation}
\label{c5b}
u'(r,a)<0 \;\;\mbox{ for }\;\; r\in (0,R(a)) \;\;\mbox{ if }\;\; a\in \mathcal{N} \cup \mathcal{N}_0\,.
\end{equation}

Let us now prove that $\mathcal{P}$ and $\mathcal{N}$ are open subsets of $(0,\infty)$. We first consider $a\in\mathcal{N}$: by \eqref{c5b} there are $\varrho>R(a)$ and $\e>0$ such that $u(\varrho,a)<0$ and $u'(r,a)<-2\e$ for $r\in (0,\varrho)$. By continuous dependence, there is $\delta\in (0,a)$ such that $u(\varrho,b)<0$ and $u'(r,b)<-\e$ for $r\in (0,\varrho)$ and $b\in (a-\delta,a+\delta)$. Since $u(0,b)=b>0$, we readily deduce that, for each $b\in (a-\delta,a+\delta)$, we have $R(b)\in (0,\varrho)$ with $u'(R(b),b)<-\e<0$. Consequently, $(a-\delta,a+\delta)$ and $\mathcal{N}$ is open in $(0,\infty)$. Consider next $a\in\mathcal{P}$, $a>1$. By Lemma~\ref{lec1} and \eqref{c2}, we have $u(r,a)\ge u(r_1(a),a)\in (0,1)$ for $r\in [0,r_1(a)]$ and $E(r_1(a),a)<0$. By continuous dependence, there is $\delta>a$ such that $u(r,b)\ge u(r_1(a),a)/2>0$ for $r\in [0,r_1(a)]$, $u(r_1(a),b)\in (0,1)$, and $E(r_1(a),b)<0$ for $b\in (a-\delta,a+\delta)$. Assume now for contradiction that there is $b\in (a-\delta,a+\delta)$ such that $R(b)<\infty$. Owing to \eqref{c2}, \eqref{c3}, and the definition of $R(b)$, we obtain
$0>E(r_1(a),b)>E(R(b),b)\ge 0$ and a contradiction. Consequently, $(a-\delta,a+\delta) \subset \mathcal{P}$ and $\mathcal{P}$ is also open in $(0,\infty)$. 

We finally argue as in \cite[Lemma~15]{Kw89} to show that there is $A>0$ such that $(A,\infty)\subset\mathcal{N}$. 

Since $\mathcal{P}$ and $\mathcal{N}$ are open subsets of $(0,\infty)$, $\mathcal{N}_0=\{a_c\}$, $(0,(p+1)^{1/p}) \subset \mathcal{P}$, and $(A,\infty)\subset\mathcal{N}$, we readily conclude that $\mathcal{P}=(0,a_c)$ and $\mathcal{N}=(a_c,\infty)$.
\end{proof}

\medskip

We next study the properties of the map $a\mapsto R(a)$. An efficient tool for that purpose is the variation of $u(.,a)$ with respect to $a$ defined by
\begin{equation*}
\vartheta(r,a) := \frac{\partial u}{\partial a}(r,a)\,, \quad (r,a)\in [0,\infty)\times (0,\infty),
\end{equation*}
which solves the second order linear differential equation
\begin{equation}
\label{c9}
\begin{array}{l}
\displaystyle{\vartheta''(r,a) + \frac{d-1}{r}\ \vartheta'(r,a) + p\ u(r,a)^{p-1}\ \vartheta(r,a) = 0}\,, \quad r\in [0,\infty)\,,\\
\\
\vartheta(0,a)= 1\,, \quad \vartheta'(0,a)=0\,,
\end{array}
\end{equation}
We argue as in \cite{FQTY08,Ta03} to prove the following lemma.

\begin{lemma}\label{lec4}
If $a>a_c$, there is a unique $z(a)\in (0,R(a))$ such that
\begin{equation*}
\left\{
\begin{array}{lcl}
\vartheta(r,a)>0 & \;\mbox{ for }\; & r\in [0,z(a))\,, \vspace{0.3cm}\\
\vartheta(z(a),a)=0 & & \vspace{0.3cm} \\
\vartheta(r,a)<0 & \;\mbox{ for }\; & r\in (z(a),R(a)]\,.
\end{array}
\right.
\end{equation*}
In addition, $u(z(a),a)>1$ and the ratio $\vartheta(.,a)/u(.,a)$ is a decreasing function of $r$ on $(0,R(a))$. 
\end{lemma}

\begin{proof}[Proof of Lemma~\ref{lec4}]
Since the proof follows rather closely that of \cite{Ta03} and \cite[Lemma~2.1]{FQTY08}, we sketch it briefly for the sake of completeness. Fix $a>a_c$ and set $u=u(.,a)$ and $\vartheta=\vartheta(.,a)$ to simplify notations. We first argue as in \cite[Lemma~17]{Kw89} to show that $\vartheta$ vanishes at least once in the interval $(0,z_1(a))$, where $z_1(a)$ denotes the unique zero in $(0,R(a))$ of $u-1$. Indeed, \eqref{c1} also reads
$$
(u(r)-1)'' + \frac{d-1}{r}\ (u(r)-1)' + \frac{u(r)^p-1}{u(r)-1}\ (u(r)-1) = 0\,, \quad r\in [0,\infty)\,
$$
and $(u(r)^p-1)/(u(r)-1)\le p\ u(r)^{p-1}$ for $r\in [0,z_1(a))$. It then follows from Sturm's comparison theorem that $\vartheta$ vanishes at least once in the interval $(0,z_1(a))$. Let $z\in (0,z_1(a))$ denote the first zero of $\vartheta$.

We now aim at showing that $\vartheta$ cannot vanish once more in the interval $(z,R(a))$. To this end, we define 
\begin{equation*}
\xi(r) :=  r^{d-1}\ \left[ u'(r)\ \vartheta(r) - u(r)\ \vartheta'(r) \right] = - r^{d-1}\ u(r)^2\ \left( \frac{\vartheta}{u} \right)'(r)\,, \quad r\in [0,R(a))\,,
\end{equation*}
which encodes the monotonicity of $\vartheta/u$. It follows from \eqref{c1} and \eqref{c9} that 
\begin{equation}
\label{c12}
\xi'(r) = r^{d-1}\ \left( (p-1)\ u^p(r) + 1 \right)\ \vartheta(r)\,, \quad r\in [0,R(a))\,.
\end{equation}
Clearly, $\xi'(r)>0$ for $r\in (0,z)$ and $\xi(0)=0$, so that $\xi(r)>0$ for $r\in (0,z]$. Assume now for contradiction that there is $\varrho\in (z,R(a))$ such that 
\begin{equation*}
\xi(r)>0 \;\mbox{ for }\; r\in (0,\varrho) \;\;\mbox{ and }\;\; \xi(\varrho)=0\,.
\end{equation*}
Observing that $\vartheta'(z)<0$, we realize that, if $\vartheta(\varrho)\ge 0$, there is $\sigma\in (z,\varrho]$ such that $\vartheta(r)<0$ for $r\in (z,\sigma)$ and $\vartheta(\sigma)=0$. In that case, $\vartheta'(\sigma)\ge 0$ and thus $\xi(\sigma)= - \sigma^{d-1}\ u(\sigma)\ \vartheta'(\sigma) \le 0$, leading us to a contradiction. Consequently,
\begin{equation}
\label{c14}
\vartheta(\varrho)<0\,.
\end{equation}
We next introduce the functions
\begin{eqnarray*}
T(r) & := & \frac{2\ (u(r)^p - 1)}{(p-1)\ u(r)^p + 1}\ \xi(r) - \zeta(r)\,,\\
\zeta(r) & := & r^d\ \left[ u'(r)\ \vartheta'(r) + (u(r)^p-1)\ \vartheta(r) \right] + (d-2)\ r^{d-1}\ u'(r)\ \vartheta(r)\,, 
\end{eqnarray*}
for $r\in [0,R(a))$ and use \eqref{c1}, \eqref{c9}, and \eqref{c12} to obtain
\begin{eqnarray}
\zeta'(r) & = & 2\ r^{d-1}\ (u(r)^p-1)\ \vartheta(r)\,, \notag\\
T'(r) & = & 2p^2\ \frac{u(r)^{p-1}}{\left[ (p-1)\ u(r)^p+1 \right]^2}\ u'(r)\ \xi(r)\,,
\label{c18}
\end{eqnarray}
for $r\in [0,R(a))$. Integrating \eqref{c18} over $(0,\varrho)$ and using the negativity of $u'$ and the positivity of $\xi$ on this interval give
\begin{equation}
\label{c19}
\zeta(\varrho)=-T(\varrho)>0\,.
\end{equation}
Since $\xi(\varrho)=0$, we have $u(\varrho)\ \vartheta'(\varrho) = u'(\varrho)\ \vartheta(\varrho)$ and we have 
$$
\zeta(\varrho) = Q(\varrho)\ \frac{\vartheta(\varrho)}{u(\varrho)}\,,
$$
where
$$
Q(r) := r^d\ \left[ u'(r)^2 + u(r)^{p+1} - u(r) \right] + (d-2)\ r^{d-1}\ u(r)\ u'(r)\,, \quad r\in [0,R(a))\,.
$$
It then follows from \eqref{c14}, \eqref{c19}, and the positivity of $u$ that 
\begin{equation}
\label{c20}
Q(\varrho)<0\,.
\end{equation}

Finally, define
$$
P(r) := r^d\ \left( u'(r)^2 + 2\ \frac{u(r)^{p+1}}{p+1} - 2\ u(r) \right) + (d-2)\ r^{d-1}\ u(r)\ u'(r)
$$
for $r\in [0,R(a))$. On the one hand, we notice that
\begin{equation}
\label{c20a}
P(r) = Q(r) -u(r) - \frac{p-1}{p+1}\ u(r)^{p+1} <Q(r)\,, \quad r\in [0,R(a))\,.
\end{equation}
On the other hand, we deduce from \eqref{c1} and \eqref{c3} that 
$$
P'(r) = r^{d-1}\ u(r)\ \left( \frac{d-2}{d-1}\ u(r)^p - (d+2) \right)\,, \quad r\in [0,R(a))\,.
$$
At this point, we realize that we have necessarily $a>(d+2)(d-1)/(d-2)$ and that there is $s\in (0,R(a))$ such that $P'(r)>0$ if $r\in (0,s)$ and $P'(r)<0$ if $r\in (s,R(a))$. Since $P(0)=0$ and $P(R(a))>0$, we conclude that $P(\varrho)>0$ and then $Q(\varrho)>0$ by \eqref{c20a}. But this contradicts \eqref{c20}. We have thus established that $\xi$ is positive in $(0,R(a))$ from which Lemma~\ref{lec4} follows. 
\end{proof}

\medskip

We are now in a position to state and prove some properties of the map $a\mapsto R(a)$. 

\begin{proposition}\label{prc3}
The map $a\mapsto R(a)$ is a decreasing function on $(a_c,\infty)$ and there is $z_1>0$ such that  
\begin{equation}
\label{c7}
\lim_{a\searrow a_c} R(a) = R(a_c) \;\;\mbox{ and }\;\; \lim_{a\to\infty} a^{(p-1)/2}\ R(a) = z_1\,.
\end{equation}
\end{proposition} 

\begin{figure}[h!]
 \begin{minipage}[t]{.5\linewidth}
\centering\epsfig{figure=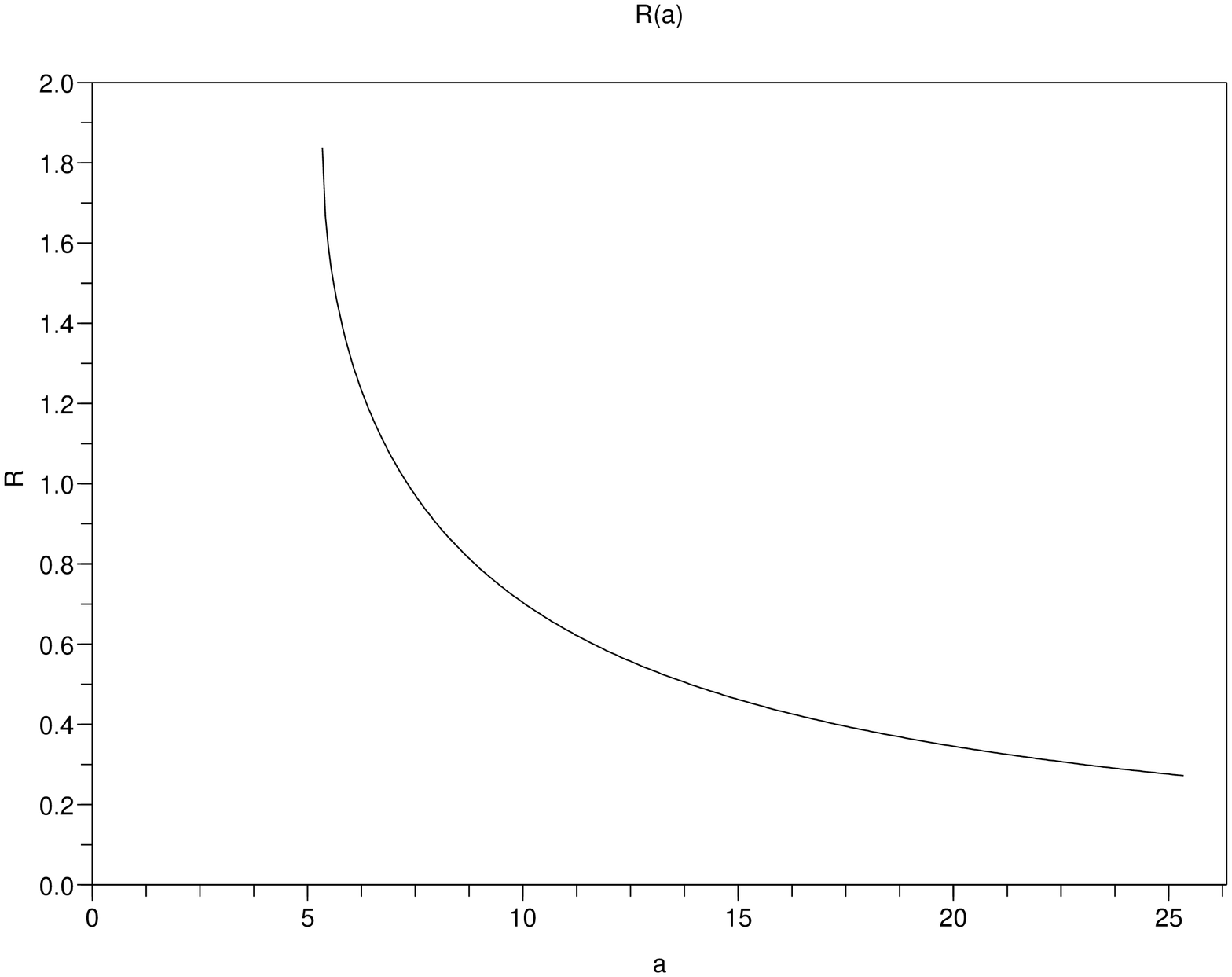,width=6cm,angle=270}
 \end{minipage}\hfill
 \begin{minipage}[t]{.5\linewidth}
\centering\epsfig{figure=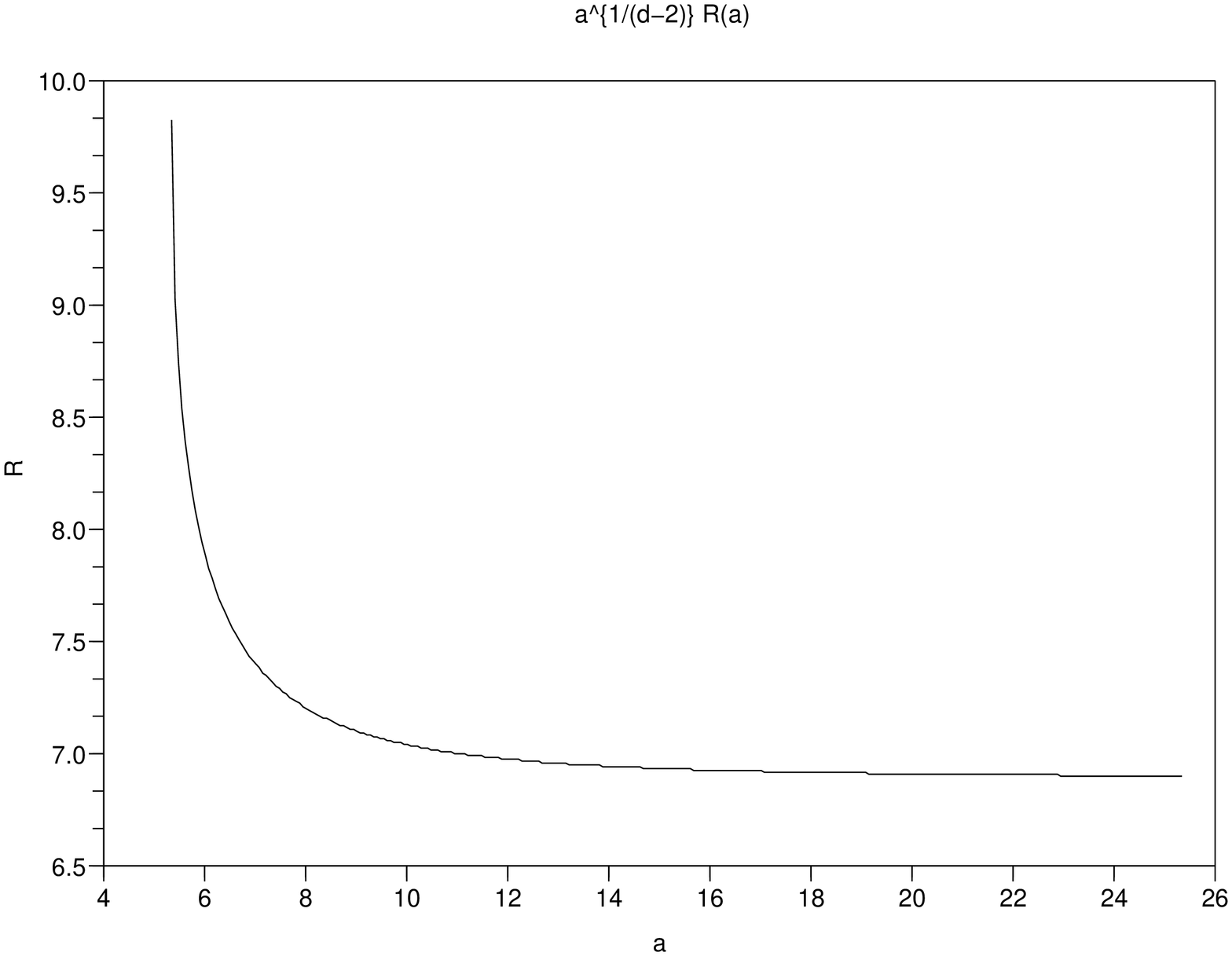,width=6cm,angle=270}
 \end{minipage}
   \caption{\small Monotonicity of the radius $R$ and $a\mapsto a^{(p-1)/2}\ R(a)$ ($d=3$).}\label{fig:prop6R}
\end{figure}

The monotonicity of $a\mapsto R(a)$ is shown in Figure~\ref{fig:prop6R}. According to numerical simulations, the function $a\mapsto a^{(p-1)/2}\ R(a)$ also seems to be a decreasing function of $a\in [a_c,\infty)$, see Figure~\ref{fig:prop6R}. 

\begin{proof}[Proof of Proposition~\ref{prc3}]
By Lemma~\ref{lec2}, $u'(R(a),a)<0$ for all $a\in (a_c,\infty)$ and the implicit function theorem warrants that $R\in\mathcal{C}^1((a_c,\infty))$ with 
$$
\frac{dR}{da}(a) = - \frac{\vartheta(R(a),a)}{u'(R(a),a)}\,.
$$
Since $\vartheta(R(a),a)<0$ by Lemma~\ref{lec4}, the previous formula implies the strict monotonicity of $a\mapsto R(a)$. We next define 
$$
R_l := \sup_{a\in (a_c,\infty)}{ R(a) }\in (0,\infty]\,.
$$
If $R_l>R(a_c)$, there is $\varrho\in (R(a_c),R_l)$ such that $u(\varrho,a_c)>0$ by Lemmata~\ref{lec1} and~\ref{lec2}. Then, there is $\delta>0$ such that $R(a)>\varrho$ for $a\in (a_c,a_c+\delta)$. It then follows from the continuous dependence of $u(.,a)$ with respect to $a$ and the monotonicity of $u(.,a)$ with respect to $r$ that
$$
0=u(R(a_c),a_c) = \lim_{a\searrow a_c} u(R(a_c),a) \ge \lim_{a\searrow a_c} u(\varrho,a) = u(\varrho,a_c)>0\,,
$$
and a contradiction. Therefore, $R_l\le R(a_c)$ is finite and we have 
$$
u(R_l,a_c)=\lim_{a\searrow a_c} u(R(a),a)= 0\,,
$$
from which we conclude that $R_l=R(a_c)$. 

Finally, define
\begin{equation}
\label{c21}
v(r,a) := \frac{1}{a}\ u\left( \frac{r}{a^{(p-1)/2}} , a \right)\,, \quad (r,a)\in [0,\infty)\times (0,\infty)\,.
\end{equation}
Owing to \eqref{c1}, $v(.,a)$ solves
\begin{equation*}
\left\{
\begin{array}{l}
\displaystyle{v''(r,a) + \frac{d-1}{r}\ v'(r,a) + |v(r,a)|^{p-1}\ v(r,a) - a^{-p} = 0}\,, \quad r\in [0,\infty)\,,\\
\\
v(0,a)= 1\,, \quad v'(0,a)=0\,,
\end{array}
\right.
\end{equation*}
In addition,
\begin{equation}
\label{c23}
v(r,a)>0 \;\;\mbox{ for }\;\; r\in \left( 0, a^{(p-1)/2}\ R(a) \right)
\end{equation} 
for $a>a_c$ by Lemma~\ref{lec2}. Since $a^{-p}\longrightarrow 0$ as $a\to\infty$, we have
\begin{equation}
\label{c23b}
\lim_{a\to\infty} \sup_{r\in [0,\varrho]}{|v(r,a)-w(r)|}=0 \;\;\mbox{ for all }\;\; \varrho>0\,,
\end{equation}
where $w$ denotes the unique solution to 
\begin{equation}
\label{c24}
\left\{
\begin{array}{l}
\displaystyle{w''(r) + \frac{d-1}{r}\ w'(r) + |w(r)|^{p-1}\ w(r) = 0}\,, \quad r\in [0,\infty)\,,\\
\\
w(0)= 1\,, \quad w'(0)=0\,.
\end{array}
\right.
\end{equation}
By \cite{Fo31}, there is $z_1>0$ such that
\begin{equation}
\label{c25}
w(r)>0 \;\mbox{ and }\; w'(r)<0\;\mbox{ for }\; r\in [0,z_1)\,, \quad w(z_1)=0\,, \quad w'(z_1)<0\,.
\end{equation}
Owing to \eqref{c25}, there is $\delta>0$ such that $w(r)<0$ for $r\in (z_1,z_1+\delta)$. It then follows from \eqref{c23b} that, given $r\in (z_1,z_1+\delta)$, $v(r,a)<0$ for $a$ large enough (depending on $r$), whence $a^{(p-1)/2}\ R(a) \le r$ for $a$ large enough by \eqref{c23}. Letting $r\to z_1$ guarantees that 
$$
\limsup_{a\to\infty} a^{(p-1)/2}\ R(a) \le z_1\,.
$$
Next, if $\varrho\in (0,\gamma)$, we have $w(r)>w(\varrho)>0$ for $r\in [0,\varrho]$ and we infer from \eqref{c23b} that $v(r,a)>w(\varrho)/2>0$ for $r\in [0,\varrho]$ and $a$ large enough. Consequently, $\varrho<a^{(p-1)/2}\ R(a)$ for $a$ large enough, from which we conclude that
$$
\liminf_{a\to\infty} a^{(p-1)/2}\ R(a) \ge z_1\,.
$$
Combining the above two inequalities completes the proof of Proposition~\ref{prc3}.
\end{proof}

\medskip

The above information allow us to estimate from above and from below a specific integral of $u(.,a)$.

\begin{proposition}\label{prc5}
For $a\in [a_c,\infty)$, we define
\begin{equation*}
\mathcal{M}(a) := d\ |B(0,1)|\ \int_0^{R(a)} u(r,a)^p\ r^{d-1}\ \dd r\,.
\end{equation*}
Recalling that $w$ is the solution to \eqref{c24} and $z_1$ is its first positive zero, we have
\begin{eqnarray}
\lim_{a\to\infty} \mathcal{M}(a) & = & \mathcal{M}_c := d\ |B(0,1)|\ \int_0^{z_1} w(r)^p\ r^{d-1}\ \dd r\,, \label{c27} \\
\mathcal{M}_2 := \sup_{a\in [a_c,\infty)}{\mathcal{M}(a)} & < & \infty\,. \label{c28}
\end{eqnarray}
\end{proposition}

\begin{proof}[Proof of Proposition~\ref{prc5}]
Let $a\ge a_c$. Since $u(0,a)=a$, it follows from the monotonicity of $u(.,a)$ that
$$
\mathcal{M}(a) \le d\ |B(0,1)|\ \int_0^{R(a)} a^p\ r^{d-1}\ \dd r = |B(0,1)|\ \left( a^{(p-1)/2}\ R(a) \right)^d\,.
$$
The upper bound \eqref{c28} is then a straightforward consequence of \eqref{c7} and the above inequality. 

Next, recalling that $v(.,a)$ is defined by \eqref{c21}, we have
$$
\mathcal{M}(a) = d\ |B(0,1)|\ \int_0^{a^{(p-1)/2} R(a)} v(r,a)^p\ r^{d-1}\ \dd r\,,
$$
and we infer from \eqref{c7} and \eqref{c23b} that \eqref{c27} holds true.
\end{proof}

\section{Proof of Theorem~\ref{th:main}}

Thanks to the analysis done in the previous sections, we are now in a position to construct self-similar blowing-up solutions to \eqref{nPKS} having either finite or infinite mass.

\begin{proposition}\label{pr:as}
Given $a>0$ and $T>0$, define
$$
\varphi(r) := \lambda_d^{d/(d-2)}\ u(\mu_d r,a)^{d/(d-2)} \;\;\mbox{ for }\;\; r\in [0,\infty) \;\;\;\mbox{ if }\;\;\; a\in (0,a_c)\,, \\
$$
and
$$
\varphi(r) := 
\left\{
\begin{array}{l}
\lambda_d^{d/(d-2)}\ u(\mu_d r,a)^{d/(d-2)} \;\;\mbox{ for }\;\; r\in [0,R(a)/\mu_d] \\
\\
0 \;\;\mbox{ for }\;\; r\ge R(a)/\mu_d\,,
\end{array}
\right. 
\;\;\;\mbox{ if }\;\;\; a\in [a_c,\infty)\,.
$$
Define next $\psi$, $\Phi$, and $\Psi$ by \eqref{bu3} and \eqref{bu21}, respectively. Then the functions $(\rho,c)$ defined by \eqref{bu1} in $(0,T)\times\RR^d$ with $s(t)=[d(T-t)]^{1/d}$ is a non-negative self-similar blowing-up solution to \eqref{nPKS} with finite mass if $a\ge a_c$ and infinite mass if $a\in (0,a_c)$.
\end{proposition}

The proof of Proposition~\ref{pr:as} readily follows from the analysis performed in Sections~\ref{bussp} and~\ref{aaode}. As for Theorem~\ref{th:main}, it is a straightforward consequence of Proposition~\ref{pr:as}, the threshold values $M_c$ and $M_2$ being given by
$$
M_c := d^{1/d}\ \left( \frac{2(d-1)}{d-2} \right)^{(d-1)/2}\ \mathcal{M}_c \;\;\mbox{ and }\;\; M_2 := d^{1/d}\ \left( \frac{2(d-1)}{d-2} \right)^{(d-1)/2}\ \mathcal{M}_2\,.
$$

\section{Discussion}

\begin{figure}[h!]
 \begin{minipage}[t]{1\linewidth}
\centering\epsfig{figure=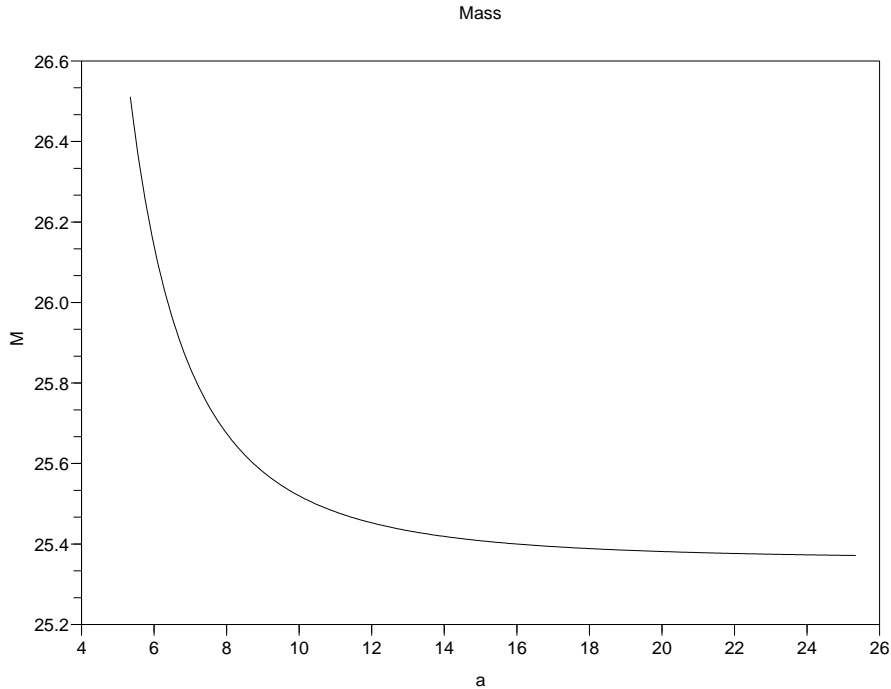,width=10cm,angle=270}
  \caption{\small Monotonicity of the mass $a\mapsto \mathcal{M}(a)$.}\label{fig:mass}
 \end{minipage} \hfill
\end{figure}

We have proved the existence of non-negative, integrable, and radially symmetric self-similar blowing-up solutions for~\eqref{nPKS}. The profile $\varphi$ of these self-similar solutions is compactly supported and non-increasing, and the mass of the corresponding self-similar solution ranges in the bounded interval $(M_c,M_2]$, the threshold mass $M_c$ corresponding to the onset of blowup found in \cite{BCL09}. Our analysis thus reveals the existence of a second threshold  value $M_2>M_c$ of the mass above which no radially symmetric and non-increasing self-similar blowing-up solution exist. The meaning of this second threshold value for the mass is yet unclear. It is worth mentioning at this point that a related situation was uncovered for the critical unstable thin-film equation
$$
\partial_t u = - \partial_x \left( u^n\ \partial^3_x u + u^{n+2}\ \partial_x u \right)\,, \quad (t,x)\in [0,\infty)\times\RR\,,
$$
in \cite{SP05} for $n\in (0,3/2)$. It is likely that, given $M \in (M_c,M_2]$, there is only a unique radially symmetric and non-increasing self-similar blowing-up solution with mass $M$ and Figure~\ref{fig:mass} provides some numerical evidence of this fact. Besides this uniqueness question, the question of stability of these blowing-up solutions is also of interest. 

\begin{figure}[h!]
 \begin{minipage}[t]{.5\linewidth}
\centering\epsfig{figure=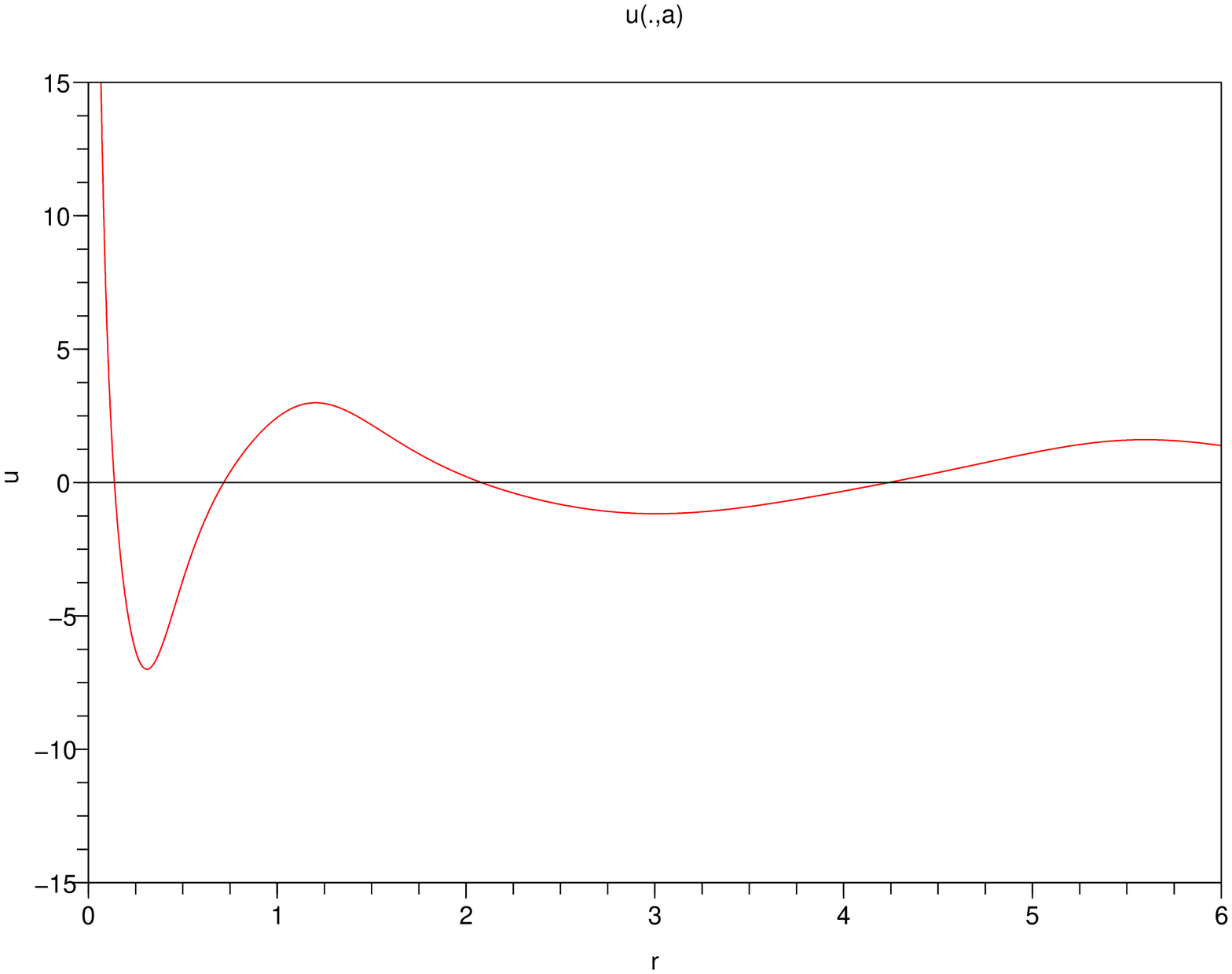,width=6cm,angle=270}
 \end{minipage}\hfill
 \begin{minipage}[t]{.5\linewidth}
\centering\epsfig{figure=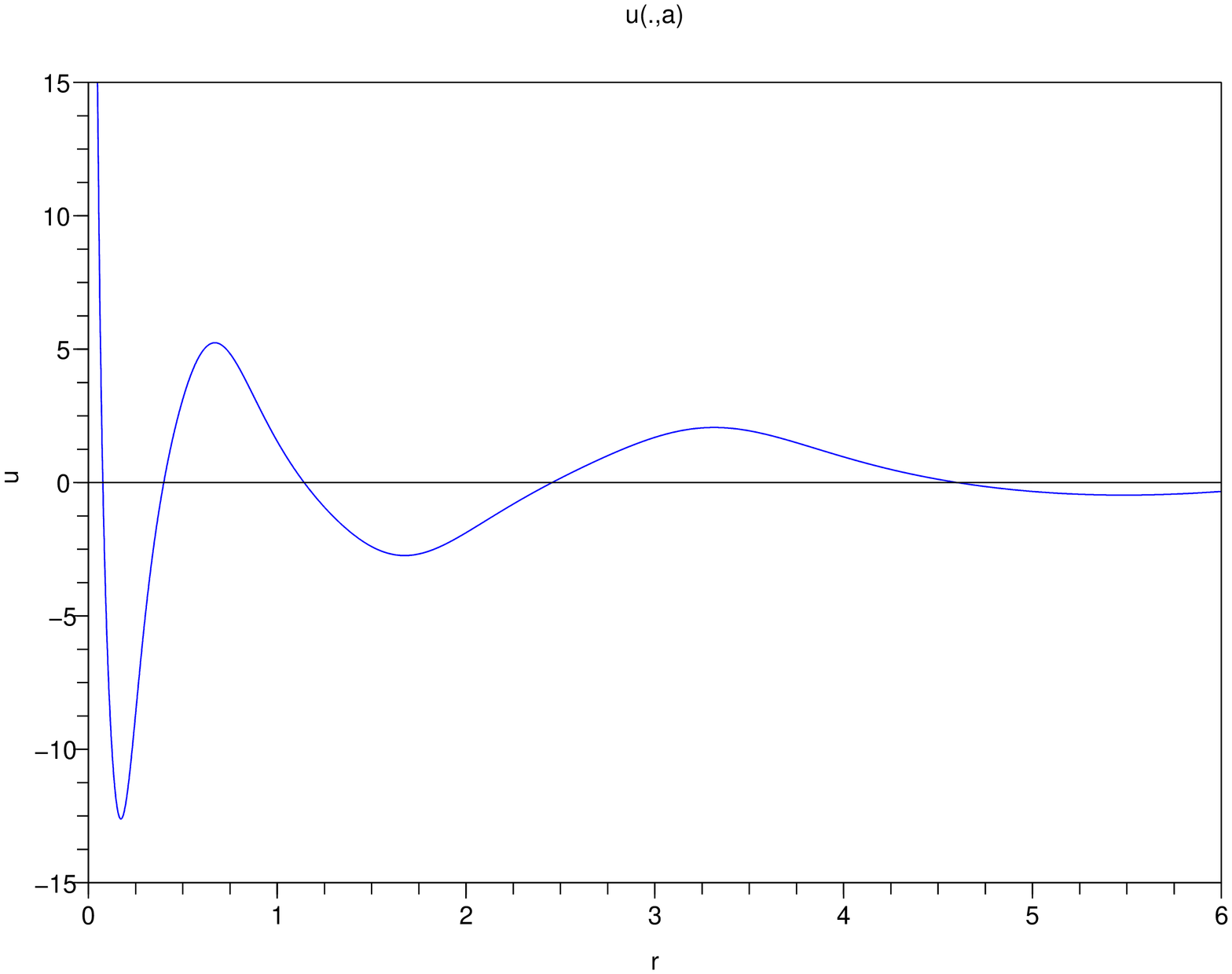,width=6cm,angle=270}
 \end{minipage}
   \caption{\small Positivity set of $u(.,a)$ with two ($a=50$, left) and three ($a=90$, right) connected components ($d=3$).}\label{fig:prc6}
\end{figure}

Another challenging question is the existence (or non-existence) of integrable profiles $\varphi$ with a non-connected positivity set as discussed in Section~\ref{bussp}. Figure~\ref{fig:prc6} provides numerical evidence that, if $a>a_c$ is large enough, $u(.,a)$ may have several zeroes and each positive ``hump'' actually corresponds to a solution of \eqref{bu9} for suitable values of $R_i$ and $R_s$. Whether the additional constraint \eqref{bu5} may be satisfied does not seem to be clear. 

\section*{Acknowledgements} We thank Jos\'e Antonio Carrillo, Jean Dolbeault, and Dejan Slep\v{c}ev for stimulating discussions. Part of this work was done while the authors enjoyed the hospitality and support of the Centro de Ciencias Pedro Pascuale de Benasque.

\bigskip\noindent{\small \copyright\, 2008 by the authors. This paper is under the Creative Commons
licence Attribution-NonCommercial-ShareAlike 2.5.}
%


\end{document}